\documentclass{article}
\usepackage{amssymb}
\usepackage{graphicx}
\usepackage{amsmath}

\newtheorem{theorem}{Theorem}

\newtheorem{conjecture}[theorem]{Conjecture}

\newtheorem{lemma}[theorem]{Lemma}

\newtheorem{question}[theorem]{Question}
\newtheorem{proposition}[theorem]{Proposition}

\newcommand{{\e}}{\ensuremath{\textbf{e}}}

\begin{document}

\title{The rank gradient from a combinatorial viewpoint}
\author{Mikl\'{o}s Ab\'{e}rt\footnote{M. Ab\'{e}rt is supported by the
grants Marie Curie IEF 235545, NSF DMS 0701105 and OTKA NK72523.}, Andrei Jaikin-Zapirain\footnote{A. Jaikin-Zapirain is supported by the Spanish Ministry of Science
and Innovation , the grant  MTM2008-06680 and the Autonomous
University of Madrid, grant CCG08-UAM/ESP-4145.},  Nikolay Nikolov}
\maketitle

\begin{abstract}
This paper investigates the asymptotic behaviour of the minimal number of
generators of finite index subgroups in residually finite groups. We analyze
three natural classes of groups: amenable groups, groups possessing an
infinite soluble normal subgroup and virtually free groups. As a tool for
the amenable case we generalize Lackenby's trichotomy theorem on finitely
presented groups.
\end{abstract}

\section{Introduction}

Let $\Gamma $ be a finitely generated group. A \emph{chain} in $\Gamma $ is a decreasing infinite sequence $\Gamma =\Gamma _{0}>\Gamma _{1}>\ldots $ of subgroups of finite index in $\Gamma $. The chain is \emph{normal} if all $\Gamma _{n}$ are normal in $\Gamma $.

For a group $\Gamma$ let $d(\Gamma)$ denote the minimal number of generators
(or rank) of $\Gamma$. For a subgroup $H\leq \Gamma $ of finite index let
\begin{equation*}
r(\Gamma ,H)=(d(H)-1)/\left\vert \Gamma :H\right\vert \text{.}
\end{equation*}%
and let the \emph{rank gradient} of $\Gamma $ with respect to the chain $%
(\Gamma _{n})$ be defined as
\begin{equation*}
\mathrm{RG}(\Gamma ,(\Gamma _{n}))=\lim_{n\rightarrow \infty }r(\Gamma
,\Gamma _{n})
\end{equation*}%
This notion has been introduced by Lackenby \cite{lack}.

In a previous paper \cite{miknik} the first and third authors investigated
the rank gradient using analytic tools, namely, the theory of cost. This tool is applicable only if the chain $(\Gamma_n)$ is normal or more generally
satisfies the Farber condition (that is, if the action of the group on
the boundary of the associated coset tree is essentially free). A
central problem discussed there is whether the rank gradient depends on the
choice of normal chain, assuming that it approximates the group. This is still unknown.

The aim of this paper is to present some results on rank gradient that can
be proved by elementary methods; that is, using just group theory and
combinatorics.

Our first three results do not assume that the chain is normal. For non-normal chains they are new and stronger than what we can show using cost.

\begin{theorem}
\label{fpamenable}Let $\Gamma $ be a finitely presented infinite amenable
group and let $(\Gamma _{n})$ be an arbitrary chain in $\Gamma $. Then $%
\mathrm{RG}(\Gamma ,(\Gamma _{n}))=0$.
\end{theorem}

This has been proved by Lackenby \cite{lack} for normal chains. 
The extension to arbitrary chains relies on a mild generalization of 
Lackenby's method plus a new ingredient, the concept of strong ergodicity for group actions discussed in Section 2. Note that the finite presentation assumption is necessary in Theorem \ref{fpamenable},
even for normal chains. Indeed, let $\Gamma =C_{2}\wr \mathbb{Z}$ be the
lamplighter group and let $\Gamma _{n}$ be the normal subgroup such that $%
\Gamma /\Gamma _{n}\cong C_{2^{n}}$. Then it is easy to see that $\Gamma _{n}
$ has a quotient equal to $C_{2}^{2^{n}}$ and so $\mathrm{RG}(\Gamma
,(\Gamma _{n}))\geq 1$. However, we do not know the answer to the following.

\begin{question}
Let $\Gamma $ be a finitely generated infinite amenable group and let $%
(\Gamma _{n})$ be a chain in $\Gamma $ with trivial intersection. Is $%
\mathrm{RG}(\Gamma ,(\Gamma _{n}))=0$?
\end{question}

In Theorem \ref{sol} we answer this question affirmatively for groups
containing an infinite soluble normal subgroup.

Theorem \ref{fpamenable} is a corollary of the following generalization of
Lackenby's theorem \cite[Theorem 1.1]{lack}.

\begin{theorem}
\label{trichotomy}Let $\Gamma $ be a finitely presented group and let $%
(\Gamma _{n})$ be an arbitrary chain in $\Gamma $. Then at least one of the
following holds: \medskip

1) the boundary action of $\Gamma $ with respect to $(\Gamma _{n})$ is
strongly ergodic; \medskip

2) $\mathrm{RG}(\Gamma ,(\Gamma _{n}))=0$; \medskip

3) there exists $n$ such that $\Gamma _{n}$ is a non-trivial amalgamated
product.
\end{theorem}

A group $\Gamma $ is a nontrivial amalgamated product $A_{1}\ast
_{A_{3}}A_{2}$ if the subgroup $A_{3}=A_{1}\cap A_{2}$ is not equal to
neither of $A_{1}$ or $A_{2}$ and has index at least $3$ in one of them.
\medskip

Compared to Lackenby's original theorem there are two new components in
Theorem \ref{trichotomy}. First, part 3) of Lackenby's theorem allows the
possibility that $\Gamma _{n}$ is an HNN extension -- we can exclude that
case.

Second, and more importantly the strong ergodicity condition replaces property ($\tau $). This is what allows us to prove Theorem \ref{fpamenable} for non-normal chains. Strong ergodicity is weaker requirement in general than property ($\tau $), it is known that for normal chains the two conditions are equivalent by
work of the first author and Elek \cite{abelek}. \medskip

From now on we will concentrate on the case when the chain has trivial
intersection.

\begin{theorem}
\label{sol} Assume that $\Gamma $ has an infinite soluble normal subgroup.
Then $\Gamma $ has rank gradient zero with respect to any chain with trivial
intersection.
\end{theorem}

Using cost, it is proved in \cite{miknik} that if $\Gamma $ has an infinite
normal amenable subgroup, then the rank gradient vanishes for normal chains
with trivial intersection. However note that when the chain is not normal Theorem \ref{sol} does not follow from the results of \cite{miknik} since the connection between cost and rank gradient established there cannot be applied. 

\bigskip

In the second half of the paper we concentrate on the case when the chain is normal, in which case stronger results are obtained. Many of the results here were known previously, however our methods are more elementary and in particular we don't use any analytic tools.

\begin{theorem}
\label{amenab} Finitely generated infinite amenable groups have rank gradient
zero with respect to any normal chain with trivial intersection.
\end{theorem}

Theorem \ref{amenab} follows from the following result of B. Weiss \cite%
{weiss}. Recall that if $N$ is a subgroup of a group $\Gamma$ then a (left) \emph {transversal} of   $N$ in   $\Gamma$ is a complete set of representatives of cosets $\{ gN:\ g\in \Gamma\}$.

\begin{theorem}
\label{weiss} Let $\Gamma $ be an amenable group generated by a finite set $%
S $ and let $(\Gamma _{n})$ be a normal chain in $\Gamma $ with trivial
intersection. Then for each $\epsilon >0$ there exists $k\in \mathbb{N}$ and
a transversal $T$ of $\Gamma_{k}$ in $\Gamma $ such that
\begin{equation*}
\left\vert TS\setminus T\right\vert <\varepsilon \left\vert T\right\vert
\text{.}
\end{equation*}
\end{theorem}

We will provide a short proof for this result of Weiss. His original proof
is a version of the Orenstein-Weiss quasitiling lemma; ours is more
algebraic and may be interesting for further applications.

Now we look at the behaviour of $r(\Gamma ,\Gamma _{n})$ over chains. It
turns out that virtually free groups can be characterized as those $\Gamma $
for which $r(\Gamma ,\Gamma _{n})$ stabilizes on normal chains $(\Gamma
_{n}) $.

\begin{theorem}
\label{stabilization} Let $\Gamma $ be a finitely generated residually
finite group. \newline
i) If $\Gamma $ is virtually free and $(\Gamma _{n})$ is a normal chain in $%
\Gamma $ with trivial intersection, then there exists $k$ such that $%
r(\Gamma ,\Gamma _{i})=r(\Gamma ,\Gamma_{k})$ ($\forall i\geq k$). \newline
ii) If $(\Gamma _{n})$ is a chain of (not necessarily normal) subgroups of $%
\Gamma $ with trivial intersection such that the sequence $r(\Gamma ,\Gamma
_{n})$ stabilizes, then $\Gamma $ is virtually free.
\end{theorem}

Next, we investigate the rank gradient of free products with amalgamation.
For free products, we obtain the following equality.

\begin{proposition}
\label{freeprod} Let $G _{1}$ and $G _{2}$ be finitely generated, residually
finite groups. Let $(N_i)$ be a normal chain in $G_1 \star G_2$ and put $%
N_{1,i}=N_i \cap G_1, \ N_{2,i}=N_i \cap G_2$.
\begin{equation*}
\mathrm{RG}(G _{1}\star G _{2},(N_i))=\mathrm{RG}(G _{1},(N_{1,i}))+\mathrm{%
RG}(G _{2},(N_{2,i}))+1.
\end{equation*}
\end{proposition}

The second theorem is very general, but it only gives an inequality.

\begin{proposition}
\label{amalgam} Let $\Gamma $ be a residually finite group generated by two
finitely generated subgroups $G _{1}$ and $G _{2}$ such that their
intersection is infinite. Then
\begin{equation*}
\mathrm{RG}(\Gamma ,(\Gamma _{n}))\leq \mathrm{RG}(G _{1},(G _{1}\cap \Gamma
_{n}))+\mathrm{RG}(G _{2},(G _{2}\cap \Gamma _{n}))
\end{equation*}%
for any normal chain $(\Gamma _{n})$ in $\Gamma $. In particular, if $G _{1}$
and $G _{2}$ have vanishing rank gradient with respect to any normal chain
then so does $\Gamma$.
\end{proposition}

Besides cost, the rank gradient is related to another important group
invariant of $\Gamma$, the first $L^{2}$ Betti number $\beta
_{1}^{(2)}(\Gamma)$. We have
\begin{equation*}
\mathrm{RG}(\Gamma ,(\Gamma _{n}))\geq \mathrm{cost}(\Gamma )-1\geq \beta
_{1}^{(2)}(\Gamma ).
\end{equation*}%
for any normal chain $(\Gamma _{n})$ in $\Gamma $ with trivial intersection.
In all the known cases so far these three numbers coincide.

There are important cases where the vanishing of the first $L^{2}$ Betti
number is known but not the cost or the rank gradient. For instance, groups
with Kazhdan property (T) have first $L^{2}$ Betti number equal to $0$ (see
\cite{bekkavalette}).

\begin{conjecture}
If $\Gamma $ has property (T) and is infinite then $\mathrm{RG}(\Gamma
,(\Gamma _{n}))=0$ for any normal chain $(\Gamma _{n})$ in $\Gamma $ with
trivial intersection.
\end{conjecture}

A group $\Gamma $ is said to be \emph{boundedly generated} if it can be
written as the product of finitely many of its cyclic subgroups $\Gamma
=\langle g_{1}\rangle \cdot \langle g_{2}\rangle \cdots \langle g_{t}\rangle
$. Examples of boundedly generated groups are arithmetic groups with the
congruence subgroup property, like $\mathrm{SL}(d,\mathbb{Z})$ ($d>2$) see
\cite{tavgen}. For many of these it follows from the results in \cite{venky}
that they have vanishing rank gradient for any normal chain. We conjecture
that in general the rank gradient of boundedly generated residually finite
groups is zero. In this direction we can show the following.

\begin{proposition}
\label{bg}If $\Gamma $ is an infinite finitely presented residually finite
boundedly generated group then the first $L^{2}$ Betti number of $\Gamma $
is zero.
\end{proposition}

The organization of the paper is as follows. In Section \ref{prelims} we
recall the Reidemeister-Schreier theorem and define coset trees and strong
ergodicity. Theorems \ref{trichotomy} and \ref{fpamenable} are proved in
Section \ref{decomp}. Section \ref{proofamenab} contains the proof of
Theorems \ref{amenab} and \ref{weiss}. Theorem \ref{sol} is proved in
Section \ref{proofsol} and Theorem \ref{stabilization}, Proposition \ref%
{freeprod} and\ \ref{amalgam} are proved in Section \ref{stabil}. Finally in
Section \ref{luck} we give a short proof of Luck's approximation theorem for
amenable groups over arbitrary fields and prove Proposition \ref{bg}.

\section{Preliminaries \label{prelims}}

First we recall the notion of Schreier graphs and the Reidemester-Schreier
theorem. Let $\Gamma $ be a group generated by a finite set $S$ and $H$ a
subgroup of finite index. Then the \emph{Schreier graph} $\Delta =\Delta
(\Gamma ,H,S)$ for $\Gamma $ relative to $H$ with respect to $S$ is an
oriented graph defined as follows:

\begin{enumerate}
\item The vertices of $\Delta$ are the left cosets of $H$ in $\Gamma$, that
is $V(\Delta)=\{gH\ |\ g\in \Gamma\}$.

\item The set of edges $E(\Delta)$ is $\{(gH,sgH)\ |\ g\in \Gamma, s\in S\}$.
\end{enumerate}

For a subset $A$ of $V(\Delta )$ we denote by $\partial A$ the set of edges
that connects $A$ and its complement $A^{c}$ in $V(\Delta )$. It is clear
that any path in the Schreier graph corresponds to a word in $\Gamma $ which
is a product of some $s_{i}$, where $s_{i}\in S\cup S^{-1}$.

Let $T$ be a left transversal for $H$ in $G$. If $g\in G$, define by $\tilde
g$ the unique $t\in T$ such that $gH=tH$. If $e=(gH, sgH)$, ($s\in S$) is an
edge of the Schreier graph $\Delta=\Delta(\Gamma,H,S)$, then we put $T(e)=({%
\widetilde{sg}})^{-1}s\tilde g $ and $T(\bar e)= \tilde g^{-1}s^{-1}{%
\widetilde{sg}} =T(e)^{-1}$. It is known that $\{T(e)\}$ generate $H$.

In this paper we will work mostly with so called Schreier transversals with
respect to $S$ for $H$ in $\Gamma $. To define them, fix a maximal tree $%
\mathcal{T}$ embedded in $\Delta $. Then for any $g\in \Gamma $ there exists
a unique path from $H$ to $gH$. Let $T$ be the set of all words
corresponding to these paths. It is clear that $T$ is a left transversal for
$H$ in $\Gamma $. We call this transversal, the \emph{Schreier transversal}
with respect to $S$ corresponding to $\mathcal{T}$. Note that $T(e)=1$ if $%
e\in E(\mathcal{T})$.

Now, let us assume that $F$ is a finitely generated free group and $S$ is a
set of its free generators. Let $H$ be a subgroup of $F$ of finite index and
put $\Delta =\Delta (F,H,S)$. Let $N$ be a normal subgroup of $F$ contained
in $H$ and generated (as a normal subgroup) by a finite set $R$. Thus, $%
F/N\cong \langle S\ |\ R\rangle $. Let $T$ be a right Schreier transversal
for $H$ in $F$ corresponding to a maximal subtree $\mathcal{T}$ of $\Delta $%
. We want to write a presentation of $H/N$ using the generators $T(e)$. Take
a relation of $F/N$, $r=s_{l}\ldots s_{1}\in R$, where $s_{i}\in S\cup S^{-1}
$ and let $t\in T$. Then $r_{t}=t^{-1}rt$ is an element of $H$ and we can
rewrite $r_{t}$ as a product of $l$ elements $T(e)$: $r_{t}=T(e_{l})\cdots
T(e_{1})$, where $e_{1}=(tH,s_{1}tH),\ldots ,e_{l}=(s_{l-1}\ldots s_{1}tH,tH)
$. Recall that $T(e)=1$ if $e\in E(\mathcal{T})$, whence we can rewrite $%
r_{t}$ as product of at most $l$ elements $T(e)^{\pm 1}$ with $e\in E(\Delta
)\setminus E(\mathcal{T})$. It is a known fact that $H/N$ has the following
presentation:
\begin{equation}
H/N\cong \left\langle \{T(e)\}_{e\in E(\Delta )\setminus E(\mathcal{T})}\ |\
\{r_{t}\}_{r\in R,t\in T}\right\rangle .  \label{presH}
\end{equation}

\bigskip

Now we define boundary actions with respect to a chain. Let $(\Gamma _{n})$
be a chain in $\Gamma $. Then the \emph{coset tree} $\mathfrak T=\mathfrak
T(\Gamma ,(\Gamma
_{n}))$ of $\Gamma $ with respect to $(\Gamma _{n})$ is defined as follows.
The vertex set of $\mathfrak T$ equals
\begin{equation*}
\mathfrak T=\left\{ g\Gamma _{n}\mid n\geq 0,g\in \Gamma \right\}
\end{equation*}%
and the edge set is defined by inclusion, that is,
\begin{equation*}
(g\Gamma _{n},h\Gamma _{m})\text{ is an edge in } \mathfrak{T}\text{ if }m=n+1\text{ and
}g\Gamma _{n}\supseteq h\Gamma _{m}
\end{equation*}%
Then $\mathfrak T$ is a tree rooted at $\Gamma $ and every vertex of level $n$ has the
same number of children, equal to the index $\left\vert \Gamma _{n}:\Gamma
_{n+1}\right\vert $. The left actions of $\Gamma $ on the coset spaces $%
\Gamma /\Gamma _{n}$ respect the tree structure and so $\Gamma $ acts on $\mathfrak T$
by automorphisms.

The boundary $\partial \mathfrak T$ of $\mathfrak T$ is defined as the set of infinite rays
starting from the root. The boundary is naturally endowed with the product
topology and product measure coming from the tree. More precisely, for $%
t=g\Gamma _{n}\in \mathfrak T$ let us define $\mathrm{Sh}(t)\subseteq \partial \mathfrak T$, the
\emph{shadow} of $t$ as
\begin{equation*}
\mathrm{Sh}(t)=\left\{ x\in \partial \mathfrak T\mid t\in x\right\}
\end{equation*}%
the set of rays going through $t$. Set the base of topology on $\partial \mathfrak T$
to be the set of shadows and set the measure of a shadow to be
\begin{equation*}
\mu (\mathrm{Sh}(t))=1/\left\vert \Gamma :\Gamma _{n}\right\vert .
\end{equation*}%
This turns $\partial \mathfrak T$ into a totally disconnected compact space with a
Borel probability measure $\mu $. The group $\Gamma $ acts ergodically on $%
\partial \mathfrak T$ by measure-preserving homeomorphisms; we call this action the
\emph{boundary action of }$\Gamma $ with respect to $(\Gamma _{n})$. See
\cite{grineksu} where these actions were first investigated in a measure
theoretic sense.

\bigskip

Let $\Gamma $ act on a probability space $(X,\mu )$ by measure preserving
maps. A sequence of subsets $A_{n}\subseteq X$ is \emph{almost invariant},
if
\begin{equation*}
\lim_{n\rightarrow \infty }\mu (A_{n}\diagdown A_{n}\gamma )=0\text{ for all
}\gamma \in \Gamma
\end{equation*}%
The sequence is trivial, if $\lim_{n\rightarrow \infty }\mu (A_{n})(1-\mu
(A_{n}))=0$. We say that the action is \emph{strongly ergodic}, if every
almost invariant sequence is trivial.

In general, spectral gap implies strong ergodicity, but not the other way
round. For a chain of subgroups $(\Gamma _{n})$ in $\Gamma $, spectral gap
is equivalent to Lubotzky's property ($\tau $), while strong ergodicity
means that large subsets of $\Gamma /\Gamma _{n}$ expand, but we do not know
what happens to small subsets.

In this paper, we will use the following two results on strongly ergodic
actions and amenability. The first is by Schmidt \cite[Theorem 2.4]{schmidt}.

\begin{theorem}
\label{amennost}Let $\Gamma $ be a countable amenable group acting on a
standard Borel probability space by measure preserving maps. Then the action
is not strongly ergodic.
\end{theorem}

The second result is from the first author and Elek \cite{abelek}.

\begin{lemma}
\label{nonstrong}Let $\Gamma $ be a group generated by a finite symmetric
set $S$ and let $(\Gamma _{n})$ be a chain in $\Gamma $ such that the
boundary action of $\Gamma $ with respect to $(\Gamma _{n})$ is not strongly
ergodic. Then for all $\varepsilon >0$ and $\alpha >\varepsilon $, for all
sufficiently large $n$ there exists a subset $A\subseteq \Gamma /\Gamma _{n}$
such that
\begin{equation*}
\left\vert \frac{\left\vert A\right\vert }{\left\vert \Gamma :\Gamma
_{n}\right\vert }-\alpha \right\vert <\varepsilon \text{ and }\left\vert
AS\setminus A\right\vert <\varepsilon \left\vert A\right\vert \text{.}
\end{equation*}
\end{lemma}

\section{Finitely presented amenable groups \label{decomp}}

First we explain the general strategy that we use to represent a finitely
presented group as an amalgamated free product. Let $H=\langle S^{H}\ |\
R^{H}\rangle $ be a finite presentation of a group $H$. We assume that any
generator from $S^{H}$ appears at least once in some relation from $R^{H}$.
Suppose $R^{H}$ is presented as a union of two subsets $R_{1}$ and $R_{2}$.
Denote by $S_{i}$ ($i=1,2$) the generators from $S$ which appears in the
words from $R_{i}$. Put $S_{3}=S_{1}\cap S_{2}$ and $R_{3}=R_{1}\cap R_{2}$.
Denote by $T_{i}$ ($i=1,2,3$) the group with the following presentation $%
\langle S_{i}\ |\ R_{i}\rangle $ and let $H_{i}$ ($i=1,2,3$) be the subgroup
of $H$ generated by $S_{i}$. It is clear that $H_{i}$ is a quotient of $%
T_{i} $.

There are two natural homomorphisms $\phi _{1}\colon T_{3}\rightarrow T_{1}$
and $\phi _{2}\colon T_{3}\rightarrow T_{2}$ which need not be injective.
Then $H$ is isomorphic to the pushout of the following diagram:
\begin{equation*}
\begin{array}{lll}
T_{3} & \rightarrow _{\phi _{1}} & T_{1} \\
\downarrow _{\phi _{2}} &  &  \\
T_{2} &  &
\end{array}%
\end{equation*}%
By the universal property of pushouts, $H_{1}\ast _{H_{3}}H_{2}$ is a
homomorphic image of $H$. Hence we have the isomorphism $H\cong H_{1}\ast
_{H_{3}}H_{2}$. Of course, in  most situations this method gives us a
trivial amalgamated free product (i.e. $H_{3}=H_{1}$ or $H_{3}=H_{2}$).

Now we give a variation of the previous construction. Suppose that $\Gamma
=\langle S\ |\ R\rangle $ and $H$ is a subgroup of finite index $\Gamma $.
Let $\Delta =\Delta (\Gamma ,H,S)$, $\mathcal{T}$ a maximal tree in $\Delta $
and $T$ the Schreier transversal for $H$ in $\Gamma $ corresponding to $%
\mathcal{T}$. Then $H$ has the following presentation (see (\ref{presH})).
\begin{equation*}
H\cong \left\langle \{T(e)\}_{e\in E(\Delta )\setminus E(\mathcal{T})}\ |\
\{r_{\tilde{g}}\}_{r\in R,\tilde{g}\in T}\right\rangle .
\end{equation*}%
We want to use the construction described in the previous paragraph. For
this we have to represent the set $R^{H}=\{r_{\tilde{g}}\}_{r\in R,\tilde{g}%
\in T}$ as a union of two subsets.

Let $A$ be a subset of $V(\Delta)$. Define $R^H(A)$ be the set of relations $%
r_{\tilde g}=T(e_1)^{\pm 1}\cdots T(e_l)^{\pm 1}$ such that for some $1\le
i\le l$ one of the end points of $e_i$ lies in $A$. Then it is clear that $%
R^H=R_1\cup R_2$ where $R_1=R^H(A)$ and $R_2= R^H(A^c)$. If all generators
from $S^H=\{T(e)\}_{e \in E(\Delta)\setminus E(\Gamma)}$ appears at least
once in some relation from $R^H$, then using the construction from the
previous paragraph we obtain the decomposition $H\cong H_1*_{H_3}H_2$. If a
generator $T(e)$ does not appear in any relation then we add it to $S_1$ if $%
e$ connects two elements from $A$, to $S_2$ if $e$ connects two elements
from $A^c$ and to $S_1$ and $S_2$ if $e$ is contained in $\partial A$. In
this case we obtain again $H\cong H_1*_{H_3}H_2$.

There is an easy description of the generating sets $S_{i}$ of $H_{i}$ ($%
i=1,2,3$). Let $X_{1}$ consists of elements $T(e)$ such that the both end
points of $e$ are in $A$, $X_{2}$ consists of elements $T(e)$ such that the
both end points of $e$ are in $A^{c}$ and $X_{3}$ consist of elements $T(e)$
such that either $e\in \partial A$ or there exists a relation $r_{\tilde{g}%
}=T(e_{1})^{\pm 1}\cdots T(e_{l})^{\pm 1}$ for which some $e_{i}\in \partial
A$ and some $e_{j}$ is equal to $e$. Then we obtain that $S_{1}=X_{1}\cup
X_{3}$, $S_{2}=X_{2}\cup X_{3}$ and $S_{3}=X_{3}$.

\bigskip

We are ready to prove the modified trichotomy result of Lackenby.

\bigskip

\noindent \textbf{Proof of Theorem \ref{trichotomy}.} We assume that $\Gamma
$ does not satisfy 1) and 2) and will show that then 3) holds. Since for $%
i\geq j$ we have
\begin{equation*}
d(\Gamma _{j})-1\geq \frac{d(\Gamma _{i})-1}{|\Gamma _{j}:\Gamma _{i}|}
\end{equation*}%
the sequence $(d(\Gamma _{i})-1)/\left\vert \Gamma :\Gamma _{i}\right\vert $
is non-increasing. Thus, using $\mathrm{RG}(\Gamma ,(\Gamma _{n}))>0$ and
changing (if needed) $\Gamma $ by $\Gamma _{j}$ we may assume that
\begin{equation*}
d(\Gamma _{i})-1\geq \frac{3|\Gamma :\Gamma _{i}|d(\Gamma )}{4}
\end{equation*}%
for all $i$.

Fix a finite presentation $\Gamma =\langle S\ |\ R\rangle $ such that $%
|S|=d(\Gamma )$. Let $L$ be the sum of the lengths of elements from $R$.
Since the boundary action of $\Gamma $ with respect to $(\Gamma _{n})$ is
not strongly ergodic, using Lemma \ref{nonstrong} there exist $i$ and $A\in
V(\Delta (\Gamma ,\Gamma _{i},S))$ such that

\begin{enumerate}
\item $\frac{|\Gamma:\Gamma_i|}4<|A|<\frac{|\Gamma:\Gamma_i|}2$ and

\item $|\partial A|<\frac 1{2(1+L^2)}|A|.$
\end{enumerate}

Put $H=\Gamma_i$ and $\Delta=\Delta(\Gamma,H,S)$. Fix a maximal tree $%
\mathcal{T}$ in $\Delta$ and let $T$ be the right Schreier transversal
corresponding to this tree.

Now we apply the construction described at the beginning of this section. We
obtain that $H$ is isomorphic to an amalgamated free product of $H_{1}\ast
_{H_{3}}H_{2}$. We will prove that $H_{3}$ has index at least $4$ in both $%
H_{1}$ and $H_{2}$. We use the previous notation, so $H_{1}$ is generated by
$S_{1}=X_{1}\cup X_{3}$, $H_{2}$ is generated by $S_{2}=X_{2}\cup X_{3}$ and
$H_{3}$ is generated by $S_{3}=X_{3}$.

Suppose that $H_{1}$ has index at most $3$ in $H_{3}$. Then $H$ is generated
by $H_{2}$ and at most one other element and so $d(H)\leq d(H_{2})+1$. It is
easy to see that
\begin{equation*}
|X_{2}|\leq |S||A^{c}|-|A^{c}|+1=(d(\Gamma )-1)|A^{c}|+1<\frac{3(d(\Gamma
)-1)|V(\Delta) |}{4}+1.
\end{equation*}%

  Ler $r=s_{l}\ldots s_{1}$ be a relation   of $\Gamma$ of length $l$. Note that there are at most $l|\partial A|$ different lifts    $r_{\tilde{g}
}=T(e_{1})^{\pm 1}\cdots T(e_{l})^{\pm 1}$ of $r$     for which some $e_{i}\in \partial
A$. And also for each such relation of $H$ we have at most $l$ generators  $T(e)$ of $H$ which are getting into 
$X_3$. Thus, if $\{l_{i}\}$ is the set of the lengths of the relations of $\Gamma$ (so $L=\sum
l_{i})$), then we have that
\begin{equation*}
|S_{3}|=|X_{3}|\leq |\partial A|(1+\sum l_{i}^{2})\leq |\partial
A|(1+L^{2})\leq \frac{|V(\Delta) |}{2}.
\end{equation*}%
Thus, we obtain that
\begin{equation*}
\frac{3d(\Gamma )|V(\Delta) |}{4}+1\leq d(H)\leq |X_{2}|+|X_{3}|+1<\frac{%
3(d(\Gamma )-1)|V(\Delta) |}{4}+2+\frac{|V(\Delta) |}{2}.
\end{equation*}

This is a contradiction when $|V(\Delta) |>8$ and therefore $|H_{1}:H_{3}|>3$.
In the same way we obtain that $|H_{2}:H_{3}|>3$. $\square $

\bigskip

\noindent \textbf{Proof of Theorem \ref{fpamenable}.} Using Theorem \ref%
{nonstrong}, the boundary action of $\Gamma $ with respect to $(\Gamma _{n})$
is not strongly ergodic. Since $\Gamma $ is amenable, it is not virtually a nontrivial amalgamated product. Using Theorem \ref{trichotomy} now, we get that 2) must hold, that is, $\mathrm{RG}(\Gamma ,(\Gamma _{n}))=0$. $\square
$

\section{Normal chains in amenable groups\label{proofamenab}}

In this section it will be sometimes convenient to work with multisets
instead of sets. If $S$ is a multiset then $|S|$ denotes the total number of
elements of $S$ counted \emph{with} repetitions.

Let $\Gamma$ be a group generated by a finite multiset $S$ and $A$ a finite
multisubset of $\Gamma$. The boundary of $A$ with respect to $S$ is the
multiset

\begin{equation*}
\partial_S(A)=\{(a,sa)\ |\ a\in A, s\in S, sa\not \in A\}.
\end{equation*}
We say that $A$ is $\epsilon$-invariant (with respect to $S$) if $%
|\partial_S(A)|\le \epsilon|S||A|$. Recall that $\Gamma$ is amenable if for
each $n$ there exists a sequence $\{A_n\}$ of finite subsets of $\Gamma$
such that $A_n$ is $a_n$-invariant and $\lim_{n\to \infty}a_n=0$. We say
that $\{A_n\}$ is a Folner sequence. It is easy to see that the definitions
of the amenability and a Folner sequence do not depend on generating
multiset $S$.

\medskip

Next we proceed with a general lemma on coverings.

\begin{lemma}
\label{lefedes}Let $G$ be a compact topological group with normalised Haar
measure $\mu $ and let $A\subseteq G$ be a measurable subset of positive
measure. For a natural number $k$ let
\begin{equation*}
\mathrm{cov}(A,k)=\max_{\substack{ X\subseteq G  \\ \left| X\right| =k}}\mu
(AX)
\end{equation*}
where
\begin{equation*}
AX=\left\{ ax\mid a\in A,x\in X\right\}
\end{equation*}
Then
\begin{equation*}
\mathrm{cov}(A,k)\geq 1-(1-\mu (A))^{k}
\end{equation*}
In particular, for $k=\left\lceil 1/\mu (A)\right\rceil $ we have
\begin{equation*}
\mathrm{cov}(A,k)>1-\frac{1}{e}
\end{equation*}
\end{lemma}

\noindent \textbf{Proof.} First, the definition makes sense, since the
maximum is always achieved by compactness. We prove the statement using
induction on $k$. For $k=1$ the statement is trivial.

Assume the lemma holds for $k-1$; this implies that there exists a subset $%
X\subseteq G$ of size $k-1$ such that $\mu (AX)\geq 1-(1-\mu (A))^{k-1}$.
Let $B=AX$ and let us define the subset
\begin{equation*}
U=\left\{ (a,g)\in G\times G\mid a\in A\text{, }ag\in B\right\}
\end{equation*}
Now $U$ is measurable in $G\times G$ and using Fubini's theorem we have
\begin{equation*}
\mu ^{2}(U)=\int_{a\in A}\mu (a^{-1}B)=\mu (A)\mu (B)
\end{equation*}
where $\mu ^{2}$ denotes the product measure on $G\times G$. Now using
Fubini's theorem from the other side gives
\begin{equation*}
\mu ^{2}(U)=\int_{g\in G}\mu (Ag\cap B)
\end{equation*}
If for all $g\in G$ we have $\mu (Ag\cap B)>\mu (A)\mu (B)$ then
\begin{equation*}
\mu (A)\mu (B)=\int_{g\in G}\mu (Ag\cap B)>\mu (A)\mu (B)
\end{equation*}
a contradiction (we used $\mu (G)=1$). So there exists $g\in G$ such that
\begin{equation*}
\mu (Ag\cap B)\leq \mu (A)\mu (B)
\end{equation*}
which implies $\mu (Ag\backslash B)\geq \mu (A)-\mu (A)\mu (B)$. Now let $%
X^{\prime }=X\cup \{g\}$. For this $X^{\prime }$ we have
\begin{equation*}
\mu (AX^{\prime })=\mu (B)+\mu (Ag\backslash B)\geq \mu (B)(1-\mu (A))+\mu
(A)\geq 1-(1-\mu (A))^{k}
\end{equation*}
using $\mu (B)\geq 1-(1-\mu (A))^{k-1}$. So the statement of the lemma holds.

Finally setting $k=\left\lceil 1/\mu (A)\right\rceil $ we have
\begin{equation*}
\mathrm{cov}(A,k)\geq 1-(1-\mu (A))^{k}\geq 1-(1-\mu (A))^{1/\mu (A)}>1-%
\frac{1}{e}
\end{equation*}
using $0<\mu (A)\leq 1$. $\square $

\bigskip

Note that for finite groups one can get a slightly better estimate using
that the intersection has integer size. What we really need here is an
absolute constant greater than $\frac{1}{2}$.

\bigskip

We will now prove Theorem \ref{weiss} in two steps. First we will show that
there exists a $c$-invariant transversal for some $c<1$, and then iterating
the first step $k$ times we will obtain $c^{k}$-invariant transversal.
\medskip

\textbf{Step 1:} Let $\delta =\frac{0.1}{1.1\cdot e}$ and let $A$ be a $\delta $%
-invariant set with respect to $S$. Since the intersection of $\Gamma _{i}$
is trivial there exists $j\in \mathbb{N}$ such that the projections $\bar{a}%
=a\Gamma _{j}$ of the elements $a\in A$ in $\overline{\Gamma }=\Gamma
/\Gamma _{j}$ are all different and $\frac{|\overline{\Gamma }|}{|A|}>10$.
Now applying Lemma \ref{lefedes}, we obtain that there exists a a subset $X$
of $\Gamma $ of size $\lceil \frac{|\overline{\Gamma }|}{|A|}\rceil $ (in particular,  $|
\overline{\Gamma }|\le |A||X|\leq |
\overline{\Gamma }|+|A|\le 1.1 |
\overline{\Gamma }|$) such
that the size of the set $\overline{AX}$ is at least $(1-\frac{1}{e})|%
\overline{\Gamma }|$. Let $B$ be a subset of the set $AX$ such that $%
\overline{B}=\overline{AX}$. So $|B|\ge (1-\frac{1}{e})|%
\overline{\Gamma }|$ and   $|AX|-|B|\leq |A||X|-(1-\frac{1}{e})|
\overline{\Gamma }| .$ 
 Thus we obtain that
\begin{equation*}
\begin{array}{lll}
|\partial _{S}(B)| & \leq  & |\partial _{S}(AX)|+|S|(|AX|-|B|)\leq
|S||A||X|(\delta +1)-(1-\frac{1}{e})|
\overline{\Gamma }|\\
&  &  \\
& \leq  & |S||\overline{\Gamma }| (1.1\cdot \delta +0.1+\frac{1}{e})\leq  
|S|\frac{|B|}{1-\frac{1}{e}}(1.1\cdot \delta +0.1+\frac{1}{e}) \\
&  &  \\
& \leq  & \frac{1.4}{e-1}|S||B|.%
\end{array}%
\end{equation*}%
We may add some $|\overline{\Gamma }|-|B|$ elements to $B$ and obtain a
tranversal $T$ for $\Gamma _{j}$ in $\Gamma $. Then
\begin{equation*}
\begin{array}{lll}
|\partial _{S}(T)| & \leq  & |\partial _{S}(B)|+|S|(|T|-|B|)\leq |S||B|(%
\frac{1.4}{e-1}-1)+|S||T| \\
&  &  \\
& \leq  & |S||T|(1-(1-\frac{1.4}{e-1})(1-\frac{1}{e}))=\frac{2.4}{e}|S||T|.%
\end{array}%
\end{equation*}%
Thus, if we put $c=\frac{2.4}{e}$ we obtain that $T$ is a $c$-invariant.
\medskip

\textbf{Step 2:} Now suppose that for some $\sigma$ there is $k\in \mathbb{N}
$ and a $\sigma$-invariant transversal $T_1$ of $\Gamma_k$ in $\Gamma$.
Using the previous step we will show that there exists $l\ge k$ and a $%
c\sigma$-transversal $T$ for $\Gamma_l$ in $\Gamma$.

If $g\in \Gamma$ denote by $\tilde g$ the unique element   from $T_1$
such that $g\Gamma_k=\tilde g \Gamma_k$. Denote by $S_1$ the multiset $\{ (%
\widetilde{st})^{-1} st\ |\ (t,st)\in \partial_S(T_1)\}$. Then $S_1$ is a
generating multiset of $\Gamma_k$. By the previous step, there is $l\ge k$
and a transversal $T_2$ for $\Gamma_k$ in $\Gamma_l$ which is $c$-invariant
with respect to $S_1$. Put $T=T_1T_2$. Let $s\in S$, $t_1\in T_1$ and $%
t_2\in T_2$. Then for the pair $(t_1t_2,st_1t_2)$ to be in $\partial_S(T)$
it is necessary that $(t_1,st_1)\in \partial_S(T_1)$ holds together with $%
(t_2, (\widetilde{st_1})^{-1} st_1t_2)\in \partial_{S_1}(T_2)$.

Hence
\begin{equation*}
|\partial_S(T)|=|\partial_{S_1}(T_2)|\le c|S_1||T_2|\le
c\sigma|S||T_1||T_2|=c\sigma|S||T|.
\end{equation*}

Iterating this process we find for any $k\in \mathbb{N}$ a transversal $T$
to some $\Gamma _{j}$ which is $c^{k}$-invariant with respect to $S$ and
Theorem \ref{weiss} is proved. $\square $\bigskip

Now Theorem \ref{amenab} follows from Theorem \ref{weiss} trivially by
noting that the Schreier set $\{(\widetilde{st})^{-1}st\ |\ (t,st)\in
\partial _{S}(T)\}$ for the transversal $T$ from Theorem  \ref{weiss} is a generating set for $%
\Gamma _{k}$ of size $|\partial _{S}(T)|\leq \epsilon |\Gamma :\Gamma _{k}||S|$.

\section{Groups with a soluble normal subgroup \label{proofsol}}

In this section we prove Theorem \ref{sol}. We start with two preliminary
results.

\begin{proposition}
\label{finitesubs} Let $\Gamma$ be a group which has a sequence of finite
normal subgroups $A_i$ such that $|A_i|\rightarrow \infty$. Then $\mathrm{RG}%
(\Gamma, (\Gamma_j))=0$ for any chain $(\Gamma_j)$ in $\Gamma$ with trivial
intersection.
\end{proposition}

\textbf{Proof.} Let $d=d(\Gamma)$. Suppose $|A_{i}|=a_{i}$ and let $H=\Gamma_j$ be a member of
the chain $(\Gamma_{i})$ of $\Gamma $, such that $A_{i}\cap \Gamma_{j}=1$.
If $|\Gamma :A_{i}H|=a$ then $|\Gamma :H|=aa_{i}$ and $A_{i}H$ can be
generated by $(d-1)a+1 \leq da$ elements. Therefore same holds for its
homomorphic image $H$, and so $r(\Gamma, \Gamma_j)\leq \frac{d}{a_{i}}$.
Since $a_{i}\rightarrow \infty $ as $i$ increases we get $\mathrm{RG}(\Gamma
,(\Gamma_{i}))=0$. $\square$ \medskip

\begin{lemma}
\label{modulegen} Let $\Gamma$ be a $d$-generated group and $N$ a $\mathbb{Z}%
\Gamma $-module generated by $t$ elements as a module over $\mathbb{Z}\Gamma$%
. Suppose that $N_{0} \leq N$ is a $\mathbb{Z}\Gamma $-submodule of index $b$%
. Then $N_{0}$ can be generated by at most $t+(2d+1)\log b$ elements as a $%
\mathbb{Z}\Gamma $-module.
\end{lemma}

We postpone the proof to the end of this section and move to \medskip

\noindent \textbf{Proof of Theorem \ref{sol}.} Suppose $(\Gamma_i)$ is a
chain with trivial intersection in a group $\Gamma$ which has an infinite
soluble normal subgroup, call it $S$. Consider the last infinite term $A$ of
the derived series of $S$. Then $A^{\prime }$ is finite and it is easy to
see that $\mathrm{RG}(\Gamma,\Gamma_i)=0$ if and only if $\mathrm{RG}%
(\Delta,(\Delta_i))=0$ where $\Delta_i= \Gamma_i A^{\prime }/A^{\prime }$
and $\Delta=\Gamma/A^{\prime }$.

Hence by considering $\Gamma/A^{\prime }$ in place of $\Gamma$ we may assume
that $A^{\prime }=1$ and $\Gamma $ has an infinite abelian normal subgroup $%
A $. \medskip

\textbf{Case 1:} Suppose that the normal closure of every element of $A$ is
finite. Then we can find a sequence $A_{i}$ of finite subgroups $A_{i}$ of $%
A $, all normal in $\Gamma $, such that $|A_{i}|\rightarrow \infty $. By
Proposition \ref{finitesubs} we are done.

\medskip

\textbf{Case 2:} Suppose that $A$ has an element whose normal closure in $%
\Gamma $ is infinite. Without loss of generality we may assume that $A$ is a
principal, i.e., $1$-generated $\Gamma $-module. For each $j\in \mathbb{N}$
put $M_{j}=A\cap \Gamma_{j}$ and $\overline{\Gamma} _{j}=A\Gamma_{j}$.
Suppose that $|\Gamma :\overline{\Gamma} _{j}|=a_{j}$ and $|A:M_{j}|=b_{j}$,
so that $|\Gamma :\Gamma_{j}|=a_{j}b_{j}$. Let $d=d(\Gamma)$. It follows that $d(\overline{%
\Gamma} _{j})<a_{j}d$ and therefore the same is true for $\overline{\Gamma}%
_j \simeq \Gamma_{j}/M_{j}$. Let $\Gamma_j=\langle h_{1},\ldots,
h_{p}\rangle M_j$ with $p<a_{j}d$.

Since $A$ is a principal $\Gamma$-module it is $a_j$-generated as a $%
\Gamma_j $-module. Now $M_j$ is a $\Gamma_j$ invariant subgroup of $A$ of
index $b_j$. By Lemma \ref{modulegen} $M_{j}$ is generated by at most $%
a_{j}+(2d(\overline \Gamma_j)+1)\log b_{j}$ elements as a module over $%
\overline \Gamma_j$. Let $M_j=\langle n_{1},\ldots ,n_{q}\rangle^{\overline{%
\Gamma}_j}$ with $q\le a_{j}+ (2da_j+1)\log b_{j}$.

We claim that $\langle h_{1},\ldots ,h_{p},n_{1},\ldots
,n_{q}\rangle=\Gamma_j$. Indeed, $A$ acts trivially by conjugation on $M_{j}$
and $\langle h_{1},\ldots , h_{p}\rangle A=\Gamma _{j}$. Therefore $%
M_{j}=\langle n_{1},\ldots n_{q}\rangle ^{\overline{\Gamma} _{j}}=\langle
n_{1},\ldots , n_{q}\rangle ^{\langle h_{1}, \ldots ,h_{p}\rangle }$, while $%
\Gamma_{j}=\langle h_{1},....h_p\rangle M_{j}$. It follows that $%
d(\Gamma_{j})\leq p+q$ and hence
\begin{equation*}
\frac{d(\Gamma_{j})-1}{|\Gamma :\Gamma_{j}|}\leq \frac{p+q}{a_{j}b_{j}}<%
\frac{d}{b_{j}}+\frac{1+(2d+1)\log b_{j}}{b_{j}}.
\end{equation*}
Since $\cap _{j}\Gamma_{j}=\{1\}$ the index $b_{j}=|A:A\cap \Gamma_{j}|
\rightarrow \infty$ with $j$ and so the right hand side tends to $0$.

This completes Case 2 and the proof of Theorem \ref{sol}. $\square $ \medskip

\noindent \textbf{Proof of Lemma \ref{modulegen}.} First we shall prove the
Lemma in the special case when $t=1$, i.e. $N$ is a principal $\mathbb{Z}%
\Gamma$ module. Without loss of generality we can assume that $N=\mathbb{Z}%
\Gamma $. Let $\langle g_{1},\ldots g_{d}\rangle =\Gamma $.

Pick $f_1,\ldots, f_s \in N$ whose images in $N/N_0$ generate it as an
abelian group and such that $s \leq \log b$. Let $K=\langle f_1,\ldots, f_s
\rangle \leq N$ and let $K_0=K\cap N_0$. It follows that the abelian group $%
K_0$ can be generated by at most $s$ elements, say $c_1,\ldots ,c_s \in N$.
Since $N=K+N_0$ there exists $e \in K$ such that $1 -e \in N_0$. Moreover,
for each pair of indices $(i,j)$ with $1 \leq i \leq d$, $1\leq j \leq s$
and each $\epsilon \in \{\pm 1\}$ there exists $e_{i,j, \epsilon} \in K$
such that
\begin{equation}  \label{gene}
f_j\cdot g_i^\epsilon - e_{i,j,\epsilon} \in N_0.
\end{equation}

Let $M$ be the $\mathbb{Z}\Gamma $-submodule of $N$ generated by the $%
1+(2d+1)s$ elements $1-e, \{c_{k}\}$, and $\{f_{j}\cdot g_{i}^{\epsilon
}-e_{i,j,\epsilon }\}$. We claim that $M=N_{0}$.

It is clear that $M\leq N_{0}$. For the opposite direction we first show
that $N=K+M$. It is enough to prove that each $g\in \Gamma \subset \mathbb{Z}%
\Gamma =N$ is in $K+M$. When $g=1$ we have $1=e+(1-e)$ with $e\in K$ and $%
1-e\in M$. Observe that the elements (\ref{gene}) of $M$ give $K\cdot
g_{i}^{\epsilon }\subseteq K+M$. Now $g=1\cdot g\equiv e\cdot g$ mod $M$ and
use induction on the length of the shortest expresion of $g$ as a product of
$g_{i}^{\pm 1}$.

Hence $N=K+M$ and so $N_{0}=K_{0}+M$. But $K_{0} =\langle c_1, \ldots, c_{s}
\rangle \leq M$, hence $K_{0}\leq M$ and $N_{0}=M$ as claimed. This proves
case $t=1$ of Lemma \ref{modulegen}.

\medskip

Now we can prove the general case by induction on $t$. Take a submodule $K
<M $ which is $t-1$ generated and $M/K$ is a principal $\mathbb{Z}\Gamma$%
-module. Let $b_1$ be index of $N_0 \cap K$ in $K$ and $b_2$ be the index of
$N_0K$ in $N$. Then $b_1b_2=b$ and we may assume that $K\cap N_0$ is $t-1 +
(1+2d)\log b_1$ generated and $N_0/N_0\cap K \simeq N_0K/K$ is $1+(1+2d)\log
b_2$ generated as $\mathbb{Z}\Gamma$ modules. Now since $N_0$ is an
extension of $N_0 \cap K$ by $N/(N_0 \cap K)$ it is generated by at most
\begin{equation*}
t-1 + (1+2d)\log b_1 + 1+(1+2d)\log b_2 = t+ (1+2d) \log b
\end{equation*}
elements.  Lemma \ref{modulegen} follows. $\square $

\section{Virtually free groups and free products with amalgamation \label%
{stabil}}

We start by analyzing when $r(\Gamma ,\Gamma _{i})$ stabilizes. \bigskip

\noindent \textbf{Proof of Theorem \ref{stabilization}.} Let $r_{i}=r(\Gamma
,\Gamma _{i})$.

\textbf{Part (i) }Assume that $\Gamma $ is virtually free. Then by
Bass-Serre theory $\Gamma $ has only finitely many conjugacy classes of
finite subgroups. As the groups $\Gamma_{i}$ are normal and $\cap
_{i}\Gamma_{i}=\{1\}$ there is $n\in \mathbb{N}$ such that $\Gamma_{n}$ is
torsion free. Now by the Stallings theorem \cite{stallings} every torsion
free virtually free group is free, hence $\Gamma_{n}$ is free and then
obviously $r_{i}=r_{n}$ for $i\geq n$.

\textbf{Part (ii) }Now assume that $r_{i}=r_{n}$ for all $i\geq n$. By
considering $\Gamma_{n}$ in place of $\Gamma $ we may assume that $n=0$ and $%
r_{i}=d(\Gamma )-1$ for all $i\geq 0$. We show that $\Gamma $ is free.

Let $d=d(\Gamma )$ and let $\Gamma =\langle S\rangle $, where $%
S=\{s_{1},\ldots s_{d}\}$. Consider the free group $F$ on $d$ free
generators $x_{1},\ldots x_{d}$ and the epimorphism $f:F\rightarrow \Gamma $
given by $f(x_{i})=s_{i}$. We claim that $f$ is an isomorphism. Assume not.
Let $w=x_{i_{k}}^{\epsilon _{k}}x_{i_{k-1}}^{\epsilon _{k-1}}\cdots
x_{i_{1}}^{\epsilon _{1}}$ be the shortest nontrivial word in $\ker f.$

Consider the segments $w_{j}=x_{i_{j}}^{\epsilon _{j}}x_{i_{j-1}}^{\epsilon
_{j-1}}\cdots x_{i_{1}}^{\epsilon _{1}}$ for $1\leq j\leq k$. The choice of $%
w$ gives that the $k$ elements $f(w_{j}),\ j=1,\ldots k$ of $\Gamma $ are
all different and therefore there exists an integer $m$ such that $%
f(w_{j})^{-1}f(w_{i})\not\in \Gamma_{m}$ for all $1\leq i\not=j\leq m$. Put $%
H=\Gamma_m$, and let $\Delta=\Delta(\Gamma, H, S)$. Fix the left Shcreier
transversal $T$ of $H$ in $\Gamma$ corresponding to a maximal tree $\mathcal{%
T}$. Then by the choice of $\Gamma_m$ we may write $f(w)$ as a product of
elements from $\{T(e), e\in E(\Delta)\}$, $f(w)=T(e_{i_k})^{\epsilon_k}%
\ldots T(e_{i_1})^{\epsilon_1}$, in such way that all $e_{i_j}$ are
different. There exist $j$ such that $e_{i_j}$ is not in $E(\mathcal{T})$,
whence we obtain that $T(e_{i_j})$ may be expressed in terms of other
generators $\{T(e),e\in E(\Delta)\setminus E(\mathcal{T})\}$ of $H$. It
follows that $d(\Gamma_{m})<(d-1)\left| \Gamma :\Gamma_{m}\right| +1$ and so
$r_{m}<d-1$, a contradiction. Hence $\Gamma $ is free. $\square $

\bigskip

Now we discuss the rank gradient of free products with amalgamation.\bigskip

\textbf{Proof of Proposition \ref{freeprod}.} Let $N$ be a normal subgroup
of index $n$ in $\Gamma =G_{1}\star G_{2}$ and denote $N_{j}=N\cap G_{j}$.
Suppose that $|G_{j}:N_{j}|=|G_{j}N:N|=k_{j}$.  The Bass-Serre theory gives
us the structure of $N$: it is a free product of $n/k_{1}$ copies of $N_{1}$
with $n/k_{2}$ copies of $N_{2}$ and a free group of rank $n-\frac{n}{k_{1}}-%
\frac{n}{k_{2}}+1$. (See the proof of Proposition \ref{amalgam} below where
a similar computation is given in the case of an amalgam.) By the
Grushko-Neumann theorem (see Proposition 3.7 in \cite{LS}) we have
\begin{equation*}
d(N)=\frac{n}{k_{1}}d(N_{1})+\frac{n}{k_{2}}d(N_{2})+n-\frac{n}{k_{1}}-\frac{%
n}{k_{2}}+1.
\end{equation*}%
Hence

\begin{equation*}
\frac{d(N)-1}{|\Gamma: N|}= \frac{d(N_1)}{|G_1: N_1|} + \frac{d(N_2)}{|G_2:
N_2|} +1.
\end{equation*}

Taking the group $N$ to range over the normal chain $(\Gamma _{i})$
Proposition \ref{freeprod} follows. $\square $ \bigskip

\bigskip

\textbf{Proof of Proposition \ref{amalgam}. }Denote by $A$ the intersection
of $G_{1}$ and $G_{2}$. Let $\widetilde{\Gamma }=G_{1}\star _{A}G_{2}$ and
let $\pi $ be the natural projection $\pi :\widetilde{\Gamma }\rightarrow
\Gamma $. By a slight abuse of notation we shall identify the groups $%
G_{1},G_{2}$ and $A$ with their preimages under $\pi $.

Let $N$ be a normal subgroup of index $n$ in $\Gamma$ and take $\widetilde
N= \pi^{-1} (N)$. Then $\widetilde N \cap G_j \leq \widetilde \Gamma$ is
isomorphic to $N \cap G_j \leq \Gamma$ under $\pi$. Moreover $d(N) \leq
d(\widetilde N)$ and $|\Gamma :N|= | \widetilde \Gamma : \widetilde N|$ so
it is enough to obtain an upper bound for $(d(\widetilde N)-1)/|\widetilde
\Gamma : \widetilde N|$.

We use the techniques of Bass-Serre theory, (as explained in \cite{serrebook}
for example). The group $\widetilde{\Gamma}$ acts on a tree $T$ with a
quotient an edge $E \subset T$ with vertices $X_1$ and $X_2$ such that $G_j =%
\mathrm{Stab}_{\Gamma}(X_j)$ and $A=\mathrm{Stab}_\Gamma (E)$.

Let $a=[\widetilde{N}A:\widetilde{N}]= |A : (N \cap A)|$ and $k_j=[%
\widetilde{N}G_j:\widetilde{N}]= [G_j : (G_j \cap \widetilde{N})]$ for $%
j=1,2 $. Then $\widetilde{N}$ is the fundamental group of the graph of
groups $\mathcal{G}=\widetilde{N}\backslash T$. The graph $G$ has $[%
\widetilde{\Gamma}:G_j\widetilde{N}]=n/k_j$ vertices of type $X_j$ for $%
j=1,2 $, and it has $n/a=[\widetilde{\Gamma}: A\widetilde{N}]$ number of
edges. The stabilizers of vertices of type $X_j$ in $\mathcal{G}$ are
isomorphic to $\widetilde{N} \cap G_j \simeq N \cap G_j$ ($j=1,2$), and all
the edge stabilizers are isomorphic to $\widetilde{N} \cap A \simeq N \cap A$%
.

Now by the presentation of the fundamental group of a graph of groups in b),
page 42 of \cite{serrebook} it follows that $\widetilde{N}$ is generated by
the stabilizers of the vertices of $\mathcal{G}$ together with elements $t_y$
for each edge of $\mathcal{G}$ which lies outside some chosen maximal
spanning tree of $\mathcal{G}$. The total number of vertices of $\mathcal{G}$
is $n/k_1 + n/k_2$ and the number of edges is $n/a$. It follows that

\begin{equation*}
d(\widetilde{N}) \leq \frac{n}{k_1}d(\widetilde{N} \cap G_1) + \frac{n}{k_2}%
d(\widetilde{N} \cap G_2)+ \frac{n}{a} - \frac{n}{k_1}-\frac{n}{k_2}+1.
\end{equation*}
and hence

\begin{equation*}
\frac{d(N)-1}{[\Gamma:N]} \leq \frac{d(\widetilde N)-1}{[\widetilde{\Gamma}:%
\widetilde{N}]}\leq \frac{d(N^{(1)})-1}{[\Gamma_1: N^{(1)}]} + \frac{%
d(N^{(2)})-1}{[\Gamma_2:N^{(2)}]} + \frac{1}{a} .
\end{equation*}
where for $j=1,2$ we denote $N^{(j)}= N \cap \Gamma_j \simeq \widetilde N
\cap \Gamma_j$.

Now take the subgroup $N$ to range over the normal chain $(\Gamma_i)$ of $%
\Gamma$ Then $(\Gamma_i^{(1)})$ and $(\Gamma_i^{(2)})$ are normal chains in $%
G_1$ and $G_2$ respectively and the numbers $a= |A:(A \cap \Gamma_i)|$ tends
to infinity with $i$. The Proposition follows. $\square$

\section{Luck Approximation and boundedly generated groups \label{luck}}

L\"{u}ck approximation gives a fast proof for Proposition \ref{bg}. \bigskip

\textbf{Proof of Proposition \ref{bg}}. Let $\Gamma =\langle g_{1}\rangle
\cdot \langle g_{2}\rangle \cdots \langle g_{t}\rangle $ and let $(\Gamma
_{i})$ be a normal chain in $\Gamma $ with trivial intersection. Let $%
K_{i}=\Gamma _{i}^{\prime }\Gamma _{i}^{2}$, let $G_{i}=\Gamma /\Gamma _{i}$
and let $H_{i}=\Gamma /K_{i}$ ($i\geq 0$). Since $\Gamma $ is infinite, $%
\left\vert G_{i}\right\vert $ tends to infinity with $i$.

Now $\Gamma_{i}/K_{i}$ is an elementary Abelian $2$-group of rank $r_{i}$
and the exponent $\exp(H_i)$ of $H_{i}$ is at most $2\left| G_{i}\right| $. Hence
\begin{equation*}
\left| G_{i}\right| 2^{r_{i}}=\left| H_{i}\right| \leq \exp  (H_{i})^t\leq
(2\left| G_{i}\right| )^{t}
\end{equation*}
and so $r_{i}\leq t+(t-1)\log _{2}\left| G_{i}\right| $. Let $d_{i}$ denote
the first Betti number of $\Gamma_{i}$. Then, using a theorem of L\"{u}ck
\cite{luck}, we have
\begin{equation*}
\beta _{1}^{2}(\Gamma )=\lim_{n\rightarrow \infty }\frac{d_{n}}{\left|
G_{n}\right| }\leq \lim_{n\rightarrow \infty }\frac{r_{n}}{\left|
G_{n}\right| }\leq \lim_{n\rightarrow \infty }\frac{t+(t-1)\log _{2}\left|
G_{n}\right| }{\left| G_{n}\right| }=0
\end{equation*}

The proposition holds. $\square $

\bigskip

We finish the paper by providing a simple proof of a result of Gabor Elek on
L\"{u}ck approximation over arbitrary field for amenable groups.

\begin{theorem}
\label{elek} Let $K$ be a field and $\Gamma$ a finitely generated amenable
group and let $(\Gamma_{i})$ be a normal chain with trivial intersection in $%
\Gamma$. Suppose $A\in\mathbb{M}_{n\times m}(K[\Gamma])$ is a matrix over
the group algebra $K[\Gamma]$ and let $A_{i}$ be the image of $A$ in $%
M_{n\times m}(K[\Gamma/\Gamma_{i}])$ under the quotient map $%
\pi_{i}:\Gamma\rightarrow\Gamma/\Gamma_{i}$. Then
\begin{equation*}
\lim_{i\rightarrow\infty}\frac{\dim_{K}\ker A_{i}}{|\Gamma :\Gamma_{i}|}
\end{equation*}
exists and does not depend on the choice of the chain $(\Gamma_{i})$.
\end{theorem}

It is an important question whether one can omit the amenability assumption
in this result for $K=\mathbb{F}_{p}$. The proof of the theorem uses the
Ornstein-Weiss lemma proved in \cite{OW}. Our exposition of this result is
based on a paper of Gromov (see \cite[page 336]{G}).

\begin{lemma}[Ornstein-Weiss]
Let $\Gamma $ be a finitely generated amenable group. Let $h(\Omega )$ be a
positive function defined on finite subsets $\Omega $ of $\Gamma $ such that

\begin{enumerate}
\item $h$ is subadditive, i.e.
\begin{equation*}
h(\Omega_1\cup \Omega_2)\le h(\Omega_1)+h(\Omega_2)
\end{equation*}
for all pairs of finite subsets $\Omega_1$and $\Omega_2$ in $\Gamma$.

\item $h$ is right $\Gamma$-invariant under $\Gamma$,
\begin{equation*}
h(\Omega\gamma)=h(\Omega), \text{\ for all\ } \gamma\in \Gamma.
\end{equation*}
\end{enumerate}

Then the limit
\begin{equation*}
\lim_{i\to\infty} h(\Omega_i)/|\Omega_i|
\end{equation*}
exists for every Folner sequence $\Omega_i\subset \Gamma$. Moreover, this
limit does not depend on the choice of the Folner sequence $(\Omega_i)$.
\end{lemma}

Now let us prove Theorem \ref{elek}. For any $c=\sum_{\gamma\in \Gamma}
c_\gamma\gamma\in K[\Gamma]$ let $\mathrm {supp} (\gamma)$ be the set $%
\{\gamma:c_\gamma\ne 0\}$. Let $\Omega$ be a finite subset of $\Gamma$. By $%
\mathbf{b }\in K[\Gamma]^m$ we denote the column vector $(b_1, \ldots
,b_m)^T $ with entries $b_i \in K[\Gamma]$. We denote the function $h$ as
follows
\begin{equation*}
h(\Omega)=\dim_K \{ A \mathbf{b} \in (K[\Gamma])^n : \mathrm
{supp} (b_i)\subset \Omega\}.
\end{equation*}

It is clear that $h$ is subadditive and right invariant. Thus, by the
Ornstein-Weiss lemma, there exists $H=\lim_{i\to\infty}
h(\Omega_i)/|\Omega_i|$ for every Folner sequences $\Omega_i\subset \Gamma$.
We will see now that the limit from Theorem \ref{elek} is equal to $m-H$.

Since $\dim _{K}\ker A_{i}/\left\vert \Gamma :\Gamma _{i}\right\vert $ is
bounded, in order to prove that the sequence $\dim _{K}\ker A_{i}/\left\vert
\Gamma :\Gamma _{i}\right\vert $ tends to $m-H$, it is enough to show the
limit of any Cauchy subsequence of $\dim _{K}\ker A_{i}/\left\vert \Gamma
:\Gamma _{i}\right\vert $ is $m-H$. Thus, without loss of generality we may
assume that $\lim_{i\rightarrow \infty }\dim _{K}\ker A_{i}/\left\vert
\Gamma :\Gamma _{i}\right\vert $ exists and we want to show that it is equal
to $m-H$.

Let $S$ be a generating set of $\Gamma $ containing the supports
of all the entries of the matrix $A$. By Theorem \ref{weiss}, for
any $\epsilon >0$ there exist $j$ and a transversal $T_{j}$ of
$\Gamma _{j}$ in $\Gamma $ such that $|\partial _{S}T_{j}|\leq
 {\epsilon } |T_{j}|$. Hence
\begin{equation*}
h(T_{j})\geq \dim _{K}\mathrm {Im}A_{j}\newline
\geq \dim _{K}\{\mathbf{a}=A\mathbf{b}\ |\ \mathrm {supp}(b_{i}),\mathrm {supp}%
(a_{i})\subset T_{j}\}\newline
\geq h(T_{j})-\epsilon |T_{j}|
\end{equation*}%
where $\mathbf{a}=(a_{1},\ldots ,a_{m})^{T}$ is a column vector in $K[\Gamma
]^{m}$. This implies that $\lim_{i\rightarrow \infty }\dim _{K}\ker
A_{i}/\left\vert \Gamma :\Gamma _{i}\right\vert =m-H$. $\square $

\bigskip

\noindent \textbf{Acknowledgement.} The authors thank Laci Pyber for helpful
comments on Lemma \ref{lefedes}.

\bigskip

\bigskip

\noindent M. Ab\'{e}rt\\
Alfr\'{e}d R\'{e}nyi Institute of Mathematics, Budapest, Hungary

\bigskip
\noindent A. Jaikin-Zapirain\\
Departamento de Matem\'aticas, Universidad Aut\'onoma de Madrid and\\
Instituto de Ciencias Matem\'aticas, CSIC-UAM-UC3M-UCM

\bigskip

\noindent N. Nikolov\\
Department of Mathematics, Imperial College London, SW7 2AZ, UK.
\end{document}